	\magnification=\magstep1
	\font\caps=cmcsc10
	\font\titlefont=cmbx10 scaled \magstep2
	\def\title#1{\centerline{\titlefont#1}\bigskip}
	\def\author#1{\centerline{\caps #1}\smallskip}
	\def\Math{\centerline{Department of Mathematics}}
	\def\UT{\centerline{The University of Texas at Austin}}
	\def\Austin{\centerline{Austin, Texas 78712}}
	\def\nat{\mathop{{\rm I}\kern-.2em{\rm N}}\nolimits}
	\def\IR{\mathop{{\rm I}\kern-.2em{\rm R}}\nolimits}
	\def\trivert{|\!|\!|}
	\def\num{\mathop{\#}\nolimits}
	\def\supp{\mathop{\rm supp}\nolimits}
	\def\wtilde{\widetilde}
	\def\un#1{$\underline{\hbox{#1}}$}
	\def\sss{{\scriptstyle s}}
	\def\iitem{\itemitem}
	\def\varep{\varepsilon}
	\def\myskip{\noalign{\vskip6pt}}
	\def\frac#1#2{{\textstyle{#1\over#2}}}
	\def\np{\vfill\eject}		
	\outer\def\demo #1. #2\par{\medbreak\noindent {\it#1.\enspace}
		{\rm#2}\par\ifdim\lastskip<\medskipamount\removelastskip
		\penalty55\medskip\fi}
\topinsert\vskip.5in\endinsert
\title{An Arbitrarily Distortable Banach Space}
\author{Thomas Schlumprecht}
\Math
\UT
\Austin
\vskip.3in
{\narrower\smallskip\noindent
{\bf Abstract}.
In this work we construct a ``Tsirelson like Banach space'' 
which is arbitrarily distortable.\smallskip}
\vskip.3in 

\baselineskip=18pt		

\beginsection{1. Introduction}

We consider the following notions.

\demo Definition.
Let $X$ be an infinite dimensional Banach space, and $\|\cdot\|$ its norm. 
If $|\cdot|$ is an equivalent norm on $X$ and $\lambda >1$ we say $|\cdot|$ 
is a {\it $\lambda$-distortion of\/} $X$ if for each infinite dimensional 
subspace $Y$ of $X$ we have 
$$\sup \left\{ {|y_1|\over |y_2|} : y_1,y_2 \in Y\quad \| y_1\| = \| y_2\|
=1\right\} \ge \lambda\ .$$ 
$X$ is called {\it $\lambda$-distortable\/} if there exists a 
$\lambda$-distortion on $X$. $X$ is called {\it distortable\/} if $X$ 
is $\lambda$-distortable for some $\lambda >1$, and $X$ is called 
{\it arbitrarily distortable\/} if $X$ is $\lambda$-distortable for all 
$\lambda >1$. 

\demo Remark. 
R.C.~James [2] showed that the spaces $\ell_1$ and $c_0$ are not distortable. 
Until now these are the only known spaces which are not distortable. 

From the proof of [6, Theorem 5.2, p.145] 
it  follows that each infinite dimensional 
uniform convex Banach space which does not contain a copy of $\ell_p$, 
$1<p<\infty$, has a distortable subspace. In [8] this result was 
generalized to any infinite dimensional Banach space which does not 
contain a copy of $\ell_p$, $1\le p< \infty$, or $c_0$. 

A famous open problem (the ``distortion problem'') is the question whether 
or not $\ell_p$, $1<p<\infty$, is distortable. 

In this paper we construct a Banach space $X$ which is arbitrarily distortable.
We first want to mention the following questions which are suggested by the 
existence of such a space. 

\demo Problem. 
Is every distortable Banach space arbitrarily distortable? Is, for example, 
Tsirelson's space $T$ (as presented in [5, Example~2.e.1]) arbitrarily 
distortable? 

\beginsection{2. Construction of $X$}

We first want to introduce some notations. 

The vector space of all real valued sequences $(x_n)$ whose elements are 
eventually zero is denoted by $c_{00}$; $(e_i)$ denotes the usual unit  
vector basis of $c_{00}$, {\it i.e.}, $e_i(j) =1$ if $i=j$ and $e_i 
(j) =0$ if $i\ne j$. For $x =\sum_{i=1}^\infty \alpha_i e_i \in c_{00}$ 
the set ${\supp (x) = \{ i\in \nat : \alpha_i\ne 0\}}$ is called the 
{\it support of\/} $x$. If $E$ and $F$ are two finite subsets of $\nat$ 
we write $E<F$ if $\max (E) < \min (F)$, and for $x,y\in c_{00}$ we write 
$x<y$ if $\supp (x) < \supp (y)$. For $E\subset \nat$ and $x= \sum_{i=1}^
\infty x_ie_i \in c_{00}$ we put $E(x) : = \sum_{i\in E} x_i e_i$. 

For the construction of $X$ we need a function $f: [1,\infty) \to [1,\infty)$ 
having the properties $(f_1)$ through $(f_5)$ as stated in the following 
lemma. The verification of $(f_1),(f_2)$, and $(f_3)$ are trivial while the 
verification of $(f_4)$ and $(f_5)$ are straightforward. 

\proclaim Lemma 1. 
Let $f(x) = \log_2(x+1)$, for $x\ge1$. Then $f$ has the following properties:
\smallskip
\iitem{\rm (f$_1$)} $f(1)=1$ and $f(x) <x$ for all $x>1$, 
\iitem{\rm (f$_2$)} $f$ is strictly increasing to $\infty$, 
\iitem{\rm (f$_3$)} $\lim_{x\to\infty} {f(x)\over x^q} =0$ for all $q>0$, 
\iitem{\rm (f$_4$)} the function $g(x)={x\over f(x)}$, $x\ge1$, is concave, and 
\iitem{\rm (f$_5$)} $f(x)\cdot f(y) \ge f(x\cdot y)$ for $x,y\ge1$.
\smallskip\par

For the sequel we fix a function $f$ having the properties stated in 
Lemma~1. 

On $c_{00}$ we define by induction for each $k\in \nat_0$ a norm $|\cdot|_k$. 
For $x=\sum x_n\cdot e_n \in c_{00}$ let $|x|_0 = \max_{n\in\nat} |x_n|$. 
Assuming that $|x|_k$ is defined for some $k\in\nat_0$ we put 
$$|x|_{k+1} = \max_{{\scriptstyle \ell\in\nat\atop\scriptstyle 
E_1<E_2<\cdots < E_\ell}\atop\scriptstyle E_i\subset \nat} 
{1\over f(\ell)} \sum_{i=1}^\ell |E_i(x)|_k\ .$$ 
Since $f(1) =1$ it follows that $(|x|_k)$ is increasing for any $x\in c_{00}$ 
and since $f(\ell) >1$ for all $\ell \ge2$ it follows that $|e_i|_k=1$ 
for any $i\in \nat$ and $k\in \nat_0$. 

Finally, we put for $x\in c_{00}$ 
$$\|x\| = \max_{k\in\nat} |x|_k\ .$$ 
Then $\|\cdot\|$ is a norm on $c_{00}$ and we let $X$ be the completion 
of $c_{00}$ with respect to $\|\cdot\|$. 

The following proposition states some easy facts about $X$. 

\proclaim Proposition 2. 
\vskip1pt
a)\enspace $(e_i)$ is a $1$-subsymmetric and $1$-unconditional basis of $X$; 
i.e., for any $x= \sum_{i=1}^\infty x_ie_i\in X$, any strictly increasing 
sequence $(n_i) \subset\nat$ and any $(\varep_i)_{i\in\nat} \subset\{-1,1\}$ 
it follows that 
$$\Big\| \sum_{i=1}^\infty x_ie_i \Big\| 
= \Big\| \sum_{i=1}^\infty \varep_i x_i e_{n_i}\Big\|\ .$$ 
b)\enspace For $x\in X$ it follows that 
$$\| x\| = \max \biggl\{ |x|_0\ , \sup_{\scriptstyle \ell \ge2 \atop 
\scriptstyle E_1<E_2< \cdots< E_\ell} {1\over f(\ell)} 
\sum_{i=1}^\ell \| E_i(x)\| \biggr\}$$ 
(where $|x|_0 = \sup_{n\in\nat} |x_n|$ for $x=\sum_{i=1}^\infty x_ie_i
\in X$). \smallskip\par

\demo Proof of Proposition 2.
\vskip1pt
Part (a) follows from the fact that $(e_i)$ is a $1$-unconditional 
and $1$-subsymmetric basis of the completion of $c_{00}$ with respect to 
$|\cdot|_k$ for any $k\in \nat_0$, which can be verified by induction 
for every $k\in \nat$.
\vskip1pt
Since $c_{00}$ is dense in $X$ it is enough to show the equation in (b) 
for an $x\in c_{00}$. If $\| x\| = |x|_0$ it follows for all $\ell\ge2$ and 
finite subsets $E_1,E_2,\ldots,E_\ell$ of $\nat$ with $E_1< E_2 < \cdots < 
E_\ell$ 
$${1\over f(\ell)} \sum_{i=1}^\ell \| E_i(x)\| 
 = \max_{k\ge 0} {1\over f(\ell)} \sum_{i=1}^\ell |E_i (x)|_k
 \le \max_{k\ge1} |x|_k \le \| x\|\ ,$$ 
which implies the assertion in this case.
\vskip1pt
If $\|x\| = |x|_k > |x|_{k-1} \ge |x|_0$, for some $k \ge 1$, there are 
$\ell,\ell' \in \nat$, $\ell\ge2$, finite subsets of $\nat$, $E_1,E_2,\ldots,
E_\ell$ and $E'_1,E'_2,\ldots,E'_{\ell'}$ with $E_1< E_2 < \cdots < E_\ell$ 
and $E'_1< E'_2< \cdots < E'_{\ell'}$, and a $k'\in \nat$ so that 
$$\eqalign{\| x\| & = |x|_k\cr
\myskip
&= {1\over f(\ell)} \sum_{i=1}^\ell |E_i(x)|_{k-1}\cr
\myskip
&\le {1\over f(\ell)} \sum_{i=1}^\ell \| E_i(x)\|\cr
\myskip
&\le \sup_{\scriptstyle 2\le \tilde \ell\atop\scriptstyle \wtilde E_1 < 
\wtilde E_2 < \cdots  \wtilde E_{\tilde\ell}} {1\over f(\tilde\ell)} 
\sum_{i=1}^{\tilde\ell} \| \wtilde E_i (x)\|\cr
\myskip
& = {1\over f(\ell')} \sum_{i=1}^{\ell'} \| E'_i (x)\|\cr
\myskip
&= {1\over f(\ell')} \sum_{i=1}^{\ell'} |E'_i(x)|_{k'}\cr
\myskip
&\le |x|_{k'+1} \le \| x\|\ ,\cr}$$ 
which implies the assertion. 

\demo Remark. 
\vskip1pt
a)\enspace  The equation in Proposition 2(b) determines the norm $\|\cdot\|$ 
in the following sense: If $\trivert \cdot\trivert$ is a norm on $c_{00}$
with $\trivert e_i\trivert =1$ for all $i\in\nat$ and with the property that 
$$\trivert x\trivert = \max \biggl\{ |x|_0\ ,\sup_{\scriptstyle \ell\ge2 \atop 
\scriptstyle E_1< E_2\cdots E_\ell} 
{1\over f(\ell)} \sum_{i=1}^\ell \trivert E_i(x)\trivert\biggr\}$$ 
for all $x\in c_{00}$, then it follows that $\|\cdot\|$ and 
$\trivert \cdot\trivert$ are equal.
Indeed one easily shows by induction for each $m\in\nat$ and each $x\in c_{00}$ 
with $\num \supp (x) =m$ that $\|x\| = \trivert x\trivert$.
\vskip1pt
b)\enspace The equation in Proposition 2(b) is similar to the equation which 
defines Tsirelson's space $T$ [LT, Example~2.e.1]. Recall that $T$ is generated 
by a norm $\| \cdot\|_T$ on $c_{00}$ satisfying the equation 
$$\|x\|_T = \max \biggl\{ |x|_0\ , \sup_{\scriptstyle \ell\in\nat \atop 
\scriptstyle \ell\le E_1 < \cdots  E_\ell} 
\frac12 \sum_{i=1}^\ell \| E_i(x)\|_T\biggr\}$$ 
(where $\ell \le E_1$ means that $\ell \le \min E_1$). Note that in the 
above equation the supremum is taken over all ``admissible collections'' 
$E_1 < E_2<\cdots < E_\ell$ (meaning that $\ell \le E_1$) while the norm 
on $X$ is computed by taking all collections $E_1< E_2< \cdots < E_\ell$. 
This forces the unit vectors in $T$ to be not subsymmetric, unlike in $X$. 
The admissibility condition, on the other hand, is necessary in order to 
imply that $T$ does not contain any $\ell_p$, $1\le p<\infty$, or $c_0$,
which was the purpose of its construction.

We will show that $X$ does not contain any subspace isomorphic to $\ell_p$, 
$1<p<\infty$, or $c_0$ and secondly that $X$ is distortable, which by [2] 
implies that it cannot contain a copy of $\ell_1$ either. Thus, in the 
case of $X$, the fact that $X$ does not contain a copy of $\ell_1$ is 
caused by the factor ${1\over f(\ell)}$ (replacing the constant factor 
$\frac12$ in $T$) which decreases to zero for increasing $\ell$. 

In order to state the main result we define for $\ell\in\nat$, $\ell\ge2$, 
and $x\in X$ 
$$\| x\|_\ell := \sup_{E_1< E_2 < \cdots < E_\ell} 
{1\over f(\ell)} \sum_{i=1}^\ell \| E_i (x)\|\ .$$ 
For each $\ell \in\nat$, $\|\cdot\|_\ell$ is a norm on $X$ and it follows that 
$${1\over f(\ell)} \|x\| \le \| x\|_\ell \le \| x\|\ ,\ \hbox{ for }\ 
x\in X\ .$$ 

\proclaim Theorem 3. 
For each $\ell \in \nat$, each $\varep >0$, and each infinite dimensional 
subspace $Z$ of $X$ there are $z_1,z_2\in Z$ with $\| z_1\| = \|z_2\|=1$ and 
$$\| z_1 \|_\ell \ge 1-\varep\ ,\quad \hbox{ and }\quad 
\|z_2 \|_\ell \le {1+\varep \over f(\ell)}\ .$$ 
In particular, $\|\cdot\|_\ell$ is an $f(\ell)$-distortion for 
each $\ell\in \nat$.

\demo Remark. 
Considering for $n\in \nat$ the space $T_{1/n}$ (see for example [1]) 
which is the completion of $c_{00}$ under the norm $\|\cdot\|_{(T,1/n)}$ 
satisfying the equation 
$$\| x\|_{(T,1/n)} =\max\biggl\{ |x|_0\ , 
\sup_{\ell\le E_1< E_2 \cdots E_\ell} {1\over n}\cdot 
\sum_{i=1}^\ell \| E_i (x)\|_{(T,1/n)}\biggr\}$$ 
for all $x\in c_{00}$ and putting for $x\in T_{1/n}$ 
$$\trivert x\trivert_{(T,1/n)} = \sup_{E_1< E_2 <\cdots < E_n} 
\sum_{i=1}^n \| E_i (x)\|_{1/n}$$ 
E.~Odell [5] observed that $\trivert\cdot\trivert_{(T,1/n)}$ is a $c\cdot n$ 
distortion of $T_{1/n}$ (where $c$ is a universal constant). 

In order to show Theorem 3 we will state the following three lemmas, and 
leave their proof for the next section. 

\proclaim Lemma 4. 
For $n\in\nat$ it follows that 
$$\Big\| \sum_{i=1}^n e_i\Big\| = {n\over f(n)}\ .$$ 

For the statement of the next lemma we need the following notion. If $Y$ 
is a Banach space with basis $(y_i)$ and if $1\le p\le\infty$ we say that 
$\ell_p$ is {\it finitely block represented\/} in $Y$ if for any $\varep >0$ 
and any $n\in\nat$ there is a normalized block $(z_i)_{i=1}^n$ of length 
$n$ of $(y_i)$, which is $(1+\varep)$-equivalent to the unit basis of 
$\ell_p^n$ and we call $(z_i)$ a block of $(y_i)$ if 
$z_i = \sum_{j=k_{i-1}+1}^
{k_i} \alpha_j y_j$  for  $i=1,2,\ldots$
and some $0=k_0 < k_1 <\cdots$ in $\nat_0$ and $(\alpha_j)\subset \IR$. 

\proclaim Lemma 5. 
$\ell_1$ is finitely block represented in each infinite block of $(e_i)$. 

\proclaim Lemma 6. 
Let $(y_n)$ be a block basis of $(e_i)$ with the following property: 
There is a strictly increasing sequence $(k_n)\subset\nat$, a sequence 
$(\varep_n)\subset \IR_+$ with $\lim_{n\to\infty} \varep_n=0$ and 
for each $n$ a normalized block basis $(y(n,i))_{i=1}^{k_n}$ which is 
$(1+\varep_n)$-equivalent to the $\ell_1^{k_n}$-unit basis so that 
$$y_n = {1\over k_n} \sum_{i=1}^{k_n} y(n,i)\ .$$ 
Then it follows for all $\ell \in\nat$ 
$$\lim_{n_1\to\infty}\ \lim_{n_2\to\infty}\ \ldots\ \lim_{n_\ell\to\infty} 
\Big\| \sum_{i=1}^\ell y_{n_i}\Big\| = {\ell \over  f(\ell)}\ .$$ 

\demo Proof of Theorem 3.
\vskip1pt
Let $Z$ be an infinite dimensional subspace of $X$ and $\varep >0$. By 
passing to a further subspace and by a standard perturbation argument we 
can assume that $Z$ is generated by a block of $(e_i)$ 
\vskip1pt
\noindent \un{Choice of $z_1$}:
\vskip1pt
By Lemma 5 and Lemma 6 one finds $(y_i)_{i=1}^\ell\subset Y$, with $y_1< 
y_2< \cdots < y_\ell$ so that $\|y_i\| \ge 1-\varep$, $1\le i\le\ell$, and 
so that $\|\sum_{i=1}^\ell y_i\| \le {\ell\over f(\ell)}$. Thus, choosing 
$$z_1 = \sum_{i=1}^\ell y_i\Big\slash \Big\| \sum_{i=1}^\ell y_i\Big\|$$ 
it follows that 
$$\eqalign{
\|z_1\|_\ell &\ge {1\over f(\ell)} \sum_{i=1}^\ell \|y_i\|\Big\slash 
\Big\| \sum_{i=1}^\ell y_i\Big\|  \qquad 
\left[ \eqalign{&\hbox{choose $E_i = \supp (y_i)$}\cr
&\hbox{for }\ i=1,\ldots,\ell\cr}\right] \cr
&\ge 1-\varep\ ,\cr}$$
which shows the desired property of $z_1$. 
\vskip1pt
\noindent \un{Choice of $z_2$}:
\vskip1pt
Let $n\in\nat$ so that ${4\ell\over n} \le\varep$ and choose according to 
Lemma~5  normalized elements $x_1<x_2<\cdots <x_n$ of $Z$ so that 
$(x_i)_{i=1}^n$ is $(1+\varep/2)$-equivalent to the unit basis of $\ell_1^n$ 
and put 
$$z_2 = \sum_{i=1}^n x_i\Big\slash \Big\| \sum_{i=1}^n x_i\Big\|\ .$$
Now let $E_1,\ldots,E_\ell$ be finite subsets of $\nat$ so that $E_1< E_2< 
\cdots < E_\ell$ and so that 
$$\| z_2\|_\ell = {1\over f(\ell)} \sum_{i=1}^\ell \| E_i(z_2)\|\ .$$ 

We can assume that $E_i$ is an interval in $\nat$ for each $i\le \ell$. For 
each $i\in\nat$ there are at most two elements $j_1,j_2\in\{1,\ldots,n\}$ 
so that $\supp(x_{j_{\sss}})\cap E_i\ne\emptyset$ 
and $\supp (x_{j_{\sss}})\setminus 
E_i\ne\emptyset$, $s=1,2$. Putting for $i=1,2,\ldots,\ell$ 
$$\wtilde E_i := \cup \bigl\{ \supp (x_j) : j\le n\ \hbox{ and }\ 
\supp (x_j) \subset E_i\bigr\}$$ 
it follows that $\| E_i (z_2)\| \le \| \wtilde E_i (z_2)\| + {2\over n}$, 
and, thus, from the fact that $(\wtilde E_i(z_2):i=1,2,\ldots,\ell)$ is a 
block of a sequence which is $(1+\varep/2)$-equivalent to the $\ell_1^n$ 
unit basis, it follows that 
$$\| z_2\|_\ell  \le {\ell \over 2n} + {1\over f(\ell)} 
\sum_{i=1}^\ell \| \wtilde E_i (z_2)\| 
\le {\varep\over2} + {1+\varep/2 \over f(\ell)} \Big\| \sum_{i=1}^\ell 
\wtilde E_i (z_2)\Big\|
\le \varep + {1\over f(\ell)} \ ,$$ 
which verifies the desired property of $z_2$. 

\beginsection{3. Proof of Lemmas 4, 5 and 6}

\demo Proof of Lemma 4.
\vskip1pt
By induction we show for each $n\in\nat$ that $\|\sum_{i=1}^n e_i\| = 
{n\over f(n)}$. If $n=1$ the assertion is clear. Assume that it is true for all 
$\tilde n<n$, where $n\ge 2$. Then there is an $\ell \in \nat$, $2\le \ell 
\le n$, and there are finite subsets of $\nat$, $E_1< E_2 < \cdots < E_\ell$, 
so that 
$$\eqalign{\Big\| \sum_{i=1}^n e_i\Big\|
&= {1\over f(\ell)} \sum_{j=1}^\ell \Big\| E_j\biggl( \sum_{i=1}^n e_i
\biggr)\Big\|\cr
\myskip
& = {1\over f(\ell)} \sum_{j=1}^\ell {n_i\over f(n_i)}\qquad
\hbox{[where $n_i=\num E_i$ and $\sum n_i=n$]}\cr
\myskip
&= {\ell\over f(\ell)} \sum_{j=1}^\ell {1\over\ell} \cdot {n_i\over f(n_i)}\cr
\myskip
&\le {\ell\over f(\ell)}\ {{n\over\ell} \over f({n\over \ell})}\qquad\quad
\hbox{[Property $(f_4)$ of Lemma 1]}\cr
\myskip
&= {n\over f(\ell)\cdot f({n\over \ell})}\cr
\myskip
&\le {n\over f(n)}\qquad\hskip.75truein
\hbox{[Property $(f_5)$ of Lemma 1]}\cr}$$
Since it is easy to see that $\|\sum_{i=1}^n e_i\| \ge {n\over f(n)}$, the 
assertion follows. 
\np

\demo Proof of Lemma 5.
\vskip1pt
The statement of Lemma 5 will essentially follow from the Theorem of Krivine 
([3] and [4]). It says that for each basic sequence $(y_n)$ there ia a 
$1\le p\le\infty$ so that $\ell_p$ is finitely block represented in 
$(y_i)$. Thus, we have to show that $\ell_p$, $1<p\le\infty$, is not finitely 
represented in any block basis of $(e_i)$. This follows from the fact that 
for any $1<p\le\infty$, any $n\in\nat$ and any block basis 
$(x_i)_{i=1}^n$ of $(e_i)$ we have
$$\Big\| {1\over n^{1/p}} \sum_{i=1}^n x_i\Big\| \ge {1\over n^{1/p}} 
\ {n\over f(n)} = {n^{1-1/p} \over f(n)}$$ 
and from $(f_3)$. 

\demo Proof of Lemma 6.
\vskip1pt
Let $y_n = {1\over k_n} \sum_{i=1}^{k_n} y(n,i)$, for $n\in\nat$ and 
$(y(n,i))_{i=1}^{k_n}$ $(1+\varep_m)$-equivalent to the $\ell_1^{k_n}$ 
unit basis.
\vskip1pt
For $x,\tilde x\in c_{00}$ and $m\in \nat$ with $x<e_m<\tilde x$ we will show 
that 
$$\lim_{n\to\infty} \| x+y_n+\tilde x^{(n)}\| = \| x+e_m + \tilde x\|\ ,
\leqno(*)$$
where 
$$\tilde x^{(n)} = \sum_{i=m+1}^\infty \tilde x_i \cdot e_{i+s_n} 
\qquad \biggl( \tilde x= \sum_{i=m+1}^\infty \tilde x_i e_i\biggr) $$
and $s_n \in\nat$ is chosen big enough so that $y_n < \tilde x^{(n)}$.
\vskip1pt
This would, together with Lemma 4, imply the assertion of Lemma~6. Indeed, 
for $\ell\in\nat$ it follows from $(*)$ that 
$$\eqalign{{\ell\over f(\ell)} & = \Big\| \sum_{i=1}^\ell e_i\Big\|
\qquad \hbox{[Lemma 4]}\cr
\myskip
&= \lim_{n\to\infty} \Big\| e_1 + \sum_{i=2}^\ell e_{i+n}\Big\|\qquad 
\hbox{[subsymmetry]}\cr
\myskip
&= \lim_{n_1\to\infty}\ \lim_{n\to\infty} \Big\| y_{n_1} + \sum_{i=2}^\ell
e_{i+n}\Big\|\cr
\myskip
&= \lim_{n_1\to\infty}\ \lim_{n\to\infty}\ \lim_{m\to\infty} 
\Big\| y_{n_1} +e_n +\sum_{i=3}^\ell e_{i+m}\Big\|\cr
\myskip
&= \lim_{n_1\to\infty}\ \lim_{n_2\to\infty}\ \lim_{m\to\infty} 
\Big\| y_{n_1}+y_{n_2} +\sum_{i=3}^\ell e_{i+m}\Big\|\cr
&\vdots\cr
&=\lim_{n_1\to\infty}\ \lim_{n_2\to\infty}\ \ldots\ \lim_{n_\ell\to\infty} 
\Big\| \sum_{i=1}^\ell y_{n_i}\Big\|\ .\cr}$$ 
In order to prove $(*)$ we show first the following 

\proclaim Claim. 
For $x,y\in c_{00}$, and $n\in\nat$, with $x<e_n<y$ and $\alpha,\beta \in
\IR_0^+$ it follows that 
$$\|x+\alpha e_n\| + \| \beta e_n+y\| \le \max \bigl\{ \| x+(\alpha +\beta) 
e_n\| + \|y\|,\|x\| + \| (\alpha+\beta) e_n+y\|\bigr\}\ .$$

We show by induction for all $k\in \nat_0$, all $x,y\in c_{00}$, and $n\in 
\nat$, with $\num\supp (x) + \num\supp (y) \le k$, and $x<e_n<y$ 
and all $q_1,q_2,\alpha,\beta \in \IR_0^+$ that 
$$\eqalign{&q_1\|x+\alpha t_n\| + q_2 \|\beta e_n+y\| \le\cr
&\qquad \le \max \bigl\{ q_1\| x+(\alpha +\beta)e_n\| 
+ q_2 \| y\|,\ q_1\|x\| + q_2 \|(\alpha+\beta) e_n+y\|\bigr\}\ .\cr}$$
For $k=0$ the assertion is trivial. Suppose it is true for some 
$k\ge 0$ and suppose $x,y\in c_{00}$, $x<e_n <y$ and $\num\supp (x) + 
\num\supp (y) = k+1$. We distinguish between the following cases. 

\demo Case 1. 
$\|x+\alpha e_n\| = |x+\alpha e_n|_0$ and $\|\beta e_n+y\| = |\beta e_n +y|_0$.

If $\|x+\alpha e_n\| = |x|_0$, then 
$$q_1\|x+\alpha e_n\| +q_2\|\beta e_n+y\|= q_1\|x\| +q_2\|\beta e_n+y\| \le 
q_1\|x\| + q_2 \| (\alpha+\beta) e_n+y\|\ .$$ 
If $\|\beta e_n +y\| = |y|_0$ we proceed similarly and if 
$\|x+\alpha e_n\|=\alpha$ and $\|\beta e_n+y\| =\beta$, and if w.l.o.g., 
$q_1\le q_2$, it follows that 
$$q_1\|x+\alpha e_n\| + q_2 \|\beta e_n+y\| = q_1\alpha +q_2\beta 
\le q_2(\alpha+\beta) \le q_1\|x\| + q_2 \|e_n (\alpha +\beta)+y\|\ .$$ 

\demo Case 2. 
$\| x+\alpha e_n\| \ne |x+\alpha e_n|_0$. 

Then we find $\ell\ge2$ and $E_1 < E_2 < \cdots < E_\ell$ so that 
$E_i\cap \supp (x)\ne \emptyset$ for $i=1,\ldots,\ell$ and 
$$\eqalign{&q_1\|x+\alpha e_n\| + q_2 \|\beta e_n+y\|\cr
\myskip
&\qquad = {q_1\over f(\ell)} \biggl[ \sum_{i=1}^{\ell-1} \| E_i(x)\| 
+ \|E_\ell (x+\alpha e_n)\|\biggr] + q_2 \|\beta e_n+y\| \cr
\myskip
&\qquad \le {q_1\over f(\ell)} \sum_{i=1}^{\ell-1} \|E_i(x)\| 
+ \left\{ \eqalign{ &{q_1\over f(\ell)} \| E_\ell (x) + (\alpha +\beta) e_n\| 
+ q_2\|y\|\cr
&\qquad \hbox{or}\cr
&{q_1\over f(\ell)} \|E_\ell (x)\| + q_2 \| (\alpha +\beta)e_n+y\|\cr}\right.
\cr
\myskip
&\qquad\hbox{[By the induction hypothesis]}\cr
&\qquad \le \max \bigl\{ q_1\|x+(\alpha+\beta) e_n\| + q_2\|y\|,\ 
q_1\|x\| + q_2 \|(\alpha +\beta) e_n+y\|\bigr\}\ ,\cr}$$ 
which shows the assertion in this case. 

In the case that $\|\beta e_n+y\| \ne |\beta e_n+y|_0$ we proceed like 
in Case~2.

In order to show the equation $(*)$ we first observe that for all $k\in 
\nat_0$, $|x+e_m+\tilde x|_k \le \| x+y_n +\tilde x^{(n)}\|$ (which easily 
follows by induction for each $k\in \nat$) and, thus, that 
$\liminf_{n\to\infty} \|x+y_n+\tilde x^{(n)}\| \ge \|x+e_m+\tilde x\|$. 
Since every subsequence of $(y_n)$ still satisfies the assumptions of 
Lemma~6 it is enough to show that 
$$\liminf_{n\to\infty} \|x+y_n+\tilde x^{(n)}\| \le \|x +e_m +\tilde x\|\ .$$
This inequality will be shown by induction for each $k\in \nat_0$ and 
all $x<e_m<\tilde x$ with $\num\supp (x) + \num\supp (\tilde x) \le k$. 
For $k=0$ the assertion is trivial. We assume the assertion to be true for 
some $k\ge0$ and we fix $x,\tilde x\in c_{00}$ with $x<e_m <\tilde x$ 
and $\num\supp (x) + \num\supp (\tilde x) = k+1$. 

We consider the following three cases: 

\demo Case 1. 
$\| x+y_n+\tilde x\| = |x+y_n+\tilde x\|_0$ for infinitely many $n\in\nat$. 
Since 
$$|x+y_n + \tilde x^{(n)}|_0 \le |x+e_m+\tilde x|_0\ ,\qquad n\in\nat \ ,$$
the assertion  follows. 

\demo Case 2. 
For a subsequence $(y'_n)$ of $(y_n)$ we have 
$$\|x+y'_n +\tilde x\| = {1\over f(\ell_n)} \sum_{i=1}^{\ell_n} \| E_i^{(n)} 
(x+ y'_n+ \tilde x)\|$$ 
where $\ell_n\uparrow \infty$ and $E_1^{(n)} < E_2^{(n)} <\cdots < 
E_{\ell_n}^{(n)}$ are finite subsets of $\nat$. Since $f(\ell_n)\to\infty$ 
for $n\to\infty$ it then follows that 
$$\liminf_{n\to\infty}  \|x+y_n +\tilde x^{(n)}\| = 1\le \| x+e_n 
+\tilde x\|\ .$$

Assume now that neither Case 1 nor Case 2 occurs. By passing to a 
subsequence we can assume 

\demo Case 3. 
There is  an $\ell \ge2$ so that 
$$\lim_{n\to\infty}\biggl( \| x+y_n+\tilde x^{(n)}\| - 
{1\over f(\ell)} \sum_{i=1}^\ell 
\| E_i^{(n)} (x+y_n +\tilde x^{(n)})\|\biggr) =0$$ 
where $E_1^{(n)} < \cdots < E_\ell^{(n)}$ are finite subsets of $\nat$ 
with the following properties: 
\smallskip
\iitem{a)} $\supp (x+y_n + \tilde x^{(n)}) \cap E_i^{(n)}\ne \emptyset$, 
$i\le \ell$, and 
\iitem{} $\supp (x+y_n+\tilde x^{(n)}) \subset \bigcup_{i=1}^\ell E_i^{(n)}$ 
\smallskip
\iitem{b)} The set $\supp (x) \cap E_i^{(n)}$, $i=1,\ldots,\ell$, does not 
depend on $n$ (note that $\supp (x) <\infty$), and we denote it by 
$E_i^x$. 
\smallskip
\iitem{c)} There are subsets $\widetilde E_1 < \widetilde E_2 < \cdots < 
\widetilde E_\ell$ of $\supp (\tilde x)$ and integers $r_n$ so that 
$\supp (\tilde x^{(n)})\cap E_i^{(n)} = \widetilde E_i + r_n$, for $n\in\nat$, 
(we use the convention that $\phi < E$ for any finite $E\subset \nat$), 
\smallskip
\iitem{d)} for $i\le \ell$ and $1\le j\le k_n$ we have either 
$\supp (y(n,j))\subset E_i^{(n)}$ or $\supp (y(n,j))\cap E_i^{(n)}=\emptyset$.
\smallskip

\noindent Indeed, letting for $i\le \ell$
$$\widetilde E_i^{(n)} := \left\{ 
\eqalign{&E_i^{(n)}\ \hbox{ if }\ E_i^n\cap \supp (y_n ) = \emptyset\cr
\myskip
&\supp \bigl(y(n,s)\bigr) \cup E_i^{(n)} \setminus \supp \bigl( y(n,t)\bigr)\cr
\myskip
&\qquad \hbox{where }\ s:= \min \bigl\{ \tilde s :\supp \bigl( y(n,\tilde s) 
\bigr) \cap E_i^{(n)}\ne \emptyset \bigr\}\cr
\myskip
&\qquad \hbox{and }\enspace t:= \max\bigl\{ \tilde s:\supp \bigl( y(n,\tilde t)
\bigr) \cap E_i^{(n)}\ne \emptyset\bigr\}\cr} 
\right.$$
the value 
$\sum_{i=1}^\ell \| E_i^{(n)} (x+y_n +\tilde x^{(n)})\|$ differs from 
$\sum_{i=1}^\ell \|\widetilde E_i^{(n)} (x+y_n +\tilde x^{(n)})\|$ at most 
by $2\ell/k_n$, which shows that d) can be assumed. 

\iitem{e)} For $i\le \ell$ the value 
$$q_i := \lim_{n\to\infty} {\num \{ j\le k_n,\supp (y(n,j))\subset E_i^{(n)}\}
\over k_n}$$
exists.\smallskip

Now we distinguish between the following subcases. 

\demo Case 3a. There are $\ell_1,\ell_2\in \nat$, so that $1\le \ell_1 \le 
\ell_2 -2< \ell_2 \le\ell$ and 
$$\eqalign{\| x+y_n + \tilde x^{(n)}\| 
&= {1\over f(\ell)} \biggl[ \sum_{i=1}^{\ell_1-1} \| E_i^{(n)} (x)\| 
+ \| E_{\ell_1}^{(n)} (x+y_n)\| 
+ \sum_{i=\ell_1+1}^{\ell_2-1} \|E_i^{(n)} (y_n)\|\cr
\myskip
&\qquad + \| E_{\ell_2}^{(n)} (y_n + \tilde x^{(n)})\| 
+ \sum_{i=\ell_2+1}^\ell \| E_i^{(n)} (\tilde x^{(n)})\|\biggr\}\ .\cr}$$
In this case it follows that 
$$\eqalign{\|x+y_n+\tilde x^{(n)}\|
&\le {1\over f(\ell)} \biggl[ \sum_{i=1}^{\ell_1} \| E_i^{(n)} (x)\| 
+ \sum_{i=\ell_1}^{\ell_2} \| E_i^{(n)} (y_n)\| 
+ \sum_{i=\ell_2}^\ell \| E_i^{(n)} (\tilde x^{(n)})\|\biggr]\cr
\myskip
&\le {1\over f(\ell)} \biggl[ \sum_{i=1}^{\ell_1} \| E_i^{(n)} (x)\| 
+ 1+\varep_n + \sum_{i=\ell_2}^\ell \| E_i^{(n)} (\tilde x^{(n)})\|
\biggr]\cr}$$
[By d) and  the fact that $(y(j,n))_{j=1}^{k_n}$ is 
$(1+\varep_n)$-equivalent to the $\ell_1^{k_n}$-unit basis]
$$\le \| x+e_m +\tilde x\| + \varep_n\quad ,$$
Note that 
$$\bigl[ \ell_1 + 1+(\ell -\ell_2+1)\le \ell\bigr]$$ 
which implies the assertion in this case. 

\demo Case 3b. There is an $1\le \ell_1\le \ell$ so that 
$$\|x+y_n +\tilde x^{(n)}\| = {1\over f(\ell)} 
\biggl[ \sum_{i=1}^{\ell_1-1} \| E_i^{(n)} (x)\| + \| E_{\ell_1}^{(n)} 
(x+y_n + \tilde x^{(n)})\| 
+ \sum_{i=\ell_1+1}^\ell \|E_i^{(n)} (x^{(n)})\| \biggr]\ .$$
Then the assertion can be deduced from the induction hypothesis (note, 
that by a) and the fact that $\ell\ge2$ we have that 
$\num\supp E_{\ell_1}^{(n)} (x+\tilde x^{(n)}) < \num\supp (x+
\tilde x^{(n)})$).

\demo Case 3c. 
There is an $\ell_1 <\ell$ so that 
$$\eqalign{\|x+y_n+\tilde x^{(n)}\| &= 
{1\over f(\ell)} \biggl[ \sum_{i=1}^{\ell_1-1} \| E_i^{(n)} (x)\| 
+ \| E_{\ell_1}^{(n)} (x+y_n) \| + \|E_{\ell_1+1}^{(n)} (y_n 
	+ \tilde x^{(n)})\|\cr
	\myskip
&\qquad+ \sum_{i=\ell_1+2}^\ell \| E_i^{(n)} (\tilde x^{(n)})\|\biggr]\ .\cr}$$
We can assume that $\supp (x) \ne 0$ and $\supp (\tilde x)\ne0$ (otherwise 
we are in case~3b). 
If $q_{\ell_1}$ (as defined in e)$\,$) vanishes it follows that 
$\lim_{n\to\infty} \| E_{\ell_1}^{(n)} (x+y_n)\| = \|E_{\ell_1}^x (x)\|$.
Otherwise there is a sequence $(j_n)\subset \nat$ with $\lim_{n\to\infty} 
j_n=\infty$ so that 
$$E_{\ell_1}^{(n)} (y_n) = {1\over k_n} \sum_{j=1}^{j_n} y(n,j)$$
and so that 
$$\lim_{n\to\infty} {j_n\over  k_n} = q_{\ell_1} >0\ .$$
Since the sequence $(E_{\ell_1}^{(n)}(y_n)/q_{\ell_1})_{n\in\nat}$ is 
asymptotically equal to the sequence $(\tilde y_n)$ with 
$\tilde y_n := {1\over j_n} \sum_{j=1}^{j_n} y(n,j)$ (note that $(\tilde y_n)$ 
satisfies the assumption of the lemma) we deduce from the induction hypothesis 
for some infinite $N\subset\nat$ that 
$$\eqalign{\lim_{\scriptstyle n\to\infty\atop\scriptstyle n\in N} 
\| E_{\ell_1}^{(n)} (x+y_n)\| 
& = q_{\ell_1} \lim_{n\to\infty} \Big\| E_{\ell_1}^x \left({x\over q_{\ell_1}}
\right) + \tilde y_n\Big\| \cr
\myskip
&\le q_{\ell_1} \Big\| E_{\ell_1}^x \left({x\over q_{\ell_1}}\right) + e_m
\Big\| \cr
\myskip
&= \| E_{\ell_1}^x (x) + q_{\ell_1} e_m\|\ .\cr}$$
Similarly we show for some infinite $M\subset N$, that 
$$\lim_{\scriptstyle n\to\infty\atop\scriptstyle n\in M} 
\| E_{\ell_1+1}^{(n)} (y_n + \tilde x^{(n)})\| 
\le \| q_{\ell_1+1} e_m + \widetilde E_{\ell_1+1} (\tilde x)\|\ .$$ 

From the claim at the beginning of the proof we deduce now that 
$$\eqalign{&\liminf_{n\to\infty} \| x+y_n+\tilde x^{(n)}\| \cr
\myskip
&\qquad \le {1\over f(\ell)} \biggl[ \sum_{i=1}^{\ell_1-1} 
\|E_i^x (x)\|+\| E_{\ell_1}^x (x) + q_{\ell_1} e_m\| + \| q_{\ell_1+1} e_m 
+ \widetilde E_{\ell_1+1} (\tilde x)\| 
+ \sum_{i=\ell_1+2}^\ell \|\widetilde E_i (\tilde x)\|\biggr]\cr
\myskip
&\qquad \le {1\over f(\ell)} \biggl[ \sum_{i=1}^{\ell_1-1}  \|E_i^x (x)\| 
+ \sum_{i=\ell_1+2}^\ell \|\widetilde E_i (\tilde x)\| \cr 
\myskip
&\qquad\qquad +\max \Bigl\{ \| E_{\ell_1}^x (x) + e_m\| 
+ \| \widetilde E_{\ell_1+1} (\tilde x)\|\ ,\ 
\| E_{\ell_1}^x (x)\| + \| e_m  + \widetilde E_{\ell_1+1} (\tilde x)\|\Bigr\}
\biggr]\cr
\myskip
&\qquad\qquad [q_{\ell_1} + q_{\ell_1+1}=1]\cr
\myskip
&\qquad \le \| x+e_m +\tilde x\|\ ,\cr}$$
which shows the assertion in this case and finishes the proof of the Lemma. 
\medskip

{\baselineskip=12pt
\frenchspacing
\noindent {\bf References}
\medskip

\item{[1]} Casazza, P.G. and Shura, Th.J., 
{\it Tsirelson's space}, Springer-Verlag, Berlin, LN~1363 (1989). 

\item{[2]} James, R.C., 
{\it Uniformly non-square Banach spaces}, 
Ann. of Math. {\bf80} (1964), 542--550.

\item{[3]} Krivine, J.L., 
{\it Sous espaces de dimension finie des espaces de Banach r\'eticul\'es}, 
Ann. of Math. {\bf104} (1976), 1--29.

\item{[4]} Lemberg, H., 
{\it Nouvelle d\'emonstration d'un th\'eor\`eme de J.L.~Krivine sum la finie 
repr\'esentation de $\ell_p$ dans un espace de Banach}, 
Israel J. Math. {\bf39} (1981), 341--348.

\item{[5]} Lindenstrauss, J. and Tzafriri, L., 
``Classical Banach Spaces I -- Sequence Spaces,''
Springer-Verlag, Berlin, 1979.

\item{[6]} Milman, V.D., 
{\it Geometric theory of Banach spaces, II: Geometry of the unit sphere}, 
Russian Math. Survey {\bf26} (1971), 79--163 (translated from Russian).

\item{[7]} Odell, E., personal communication. 

\item{[8]} Odell, E., Rosenthal, H. and Schlumprecht, Th., 
{\it On distorted norms in Banach spaces and the existence of 
$\ell^p$-types}, in preparation. 

}

\end